\documentclass [reqno,12pt]{amsart}
\usepackage{amsfonts}
\usepackage{amsmath}
\usepackage{amssymb}
\usepackage{graphicx}

\theoremstyle{plain}

\numberwithin{equation}{section} \textheight  22 true cm \textwidth  15 true cm  
\begin{document}
\title[Metallic maps between metallic Riemannian manifolds...]{%
METALLIC MAPS BETWEEN METALLIC RIEMANNIAN MANIFOLDS AND CONSTANCY OF CERTAIN MAPS}
\author{Mehmet Akif Akyol}
\address{Bing\"{o}l University, Faculty of Art and Science, Deparment of Mathematics, 12000, Bing\"{o}l, Turkey}
\email{mehmetakifakyol@bingol.edu.tr}
\subjclass[2000]{53C15, 53C20, 53C50}
\keywords{metallic Riemannian manifold, Golden Riemannian manifold, almost complex manifold, almost product manifold, holomorphic map.}

\begin{abstract}
In this paper, we introduce metallic maps between metallic Riemannian manifolds, provide an example and obtain certain conditions for such maps to be totally geodesic. We also give a sufficient condition for a map between metallic Riemannian manifolds to be harmonic map. Then we investigate the constancy of certain maps between metallic Riemannian manifolds and various manifolds by imposing the holomorphic-like condition. Moreover, we check the reverse case and show that some such maps are constant if there is a condition for this.
\end{abstract}

\maketitle

\section{Introduction}

In differential geometry, it is desirable to introduce and use suitable types of maps between
Riemannian manifolds. Such maps may help to compare geometric properties of manifolds.
Indeed, almost complex manifolds, almost contact
manifolds, almost product manifolds and almost para-contact metric manifolds and maps between such manifolds have
been studied extensively by many authors.

The number $\phi=(1 + \sqrt{5})/2 = 1,618$... which is a solution of the equation
$x^2-x-1=0,$ represents the golden ratio. Being inspired by the Golden
ratio, the notion of Golden manifold M was defined in \cite{CH2} by a tensor field $\Phi$ on M
satisfying $\Phi^2 = \Phi+I$. The authors studied properties of Golden manifolds and they
showed that $\Phi$ is an automorphism of the tangent bundle TM and its eigenvalues
are $\phi$ and $1-\phi$. Also see: \cite{CH3}, \cite{E}, \cite{GCS}, \cite{OCT}.

In \cite{BA}, we introduced a new map between Golden Riemannian manifolds by
imposing a holomorphic-like condition. We show that such map is a harmonic map,
and then we obtain a certain condition for such maps to be totally geodesic.

As a generalization of the golden mean was introduced in 1997 by vera W. de Spinadel
in \cite{S1}-\cite{S5} and called {\it metallic means family} or {\it metallic proportions}. More precisely,
fix two positive integers $p$ and $q.$ The positive solution of the equation $x^2-px-q=0$
is named member of the metallic means family \cite{S1}-\cite{S5}. These numbers, denoted by:
$$
\sigma_{p,q}=\frac{p+\sqrt{p^2+4q}}{2},
$$
are also called $(p,q)$-metallic numbers.
The members of the metallic means family have the property of carrying the name of a metal, like
the golden mean and its relatives: the silver mean, the bronze mean, the copper
mean and many others. The notion of metallic manifold $M$ was defined in \cite{CH1} by a tensor
field $J$ on $M$ satisfying $J^2=pJ+qI,$ where $p$ and $q$ positive integers.
The authors studied properties of metallic manifolds and they showed that $J$ is
an automorphism of the tangent bundle $TM$ and its eigenvalues are $\sigma_{p,q}$ and $p-\sigma_{p,q}$.
They also find conditions for this kind of submanifold to be also a metallic Riemannian manifold
in terms of invariance.

In this paper, we study a new map between metallic Riemannian manifolds by
imposing a holomorphic-like condition for the first time as far as we know.
We give an example, obtain certain conditions for such maps to be totally geodesic and
give a sufficient condition for a map between metallic Riemannian manifolds to be harmonic map.
Moreover, we also check the existence of such maps between metallic Riemannian manifolds and
another manifold equipped with a differentiable structure(Golden, almost product, almost complex,
almost contact, almost para contact).\\

\section{Preliminaries}

In this section, we give a brief information for almost complex manifolds,
almost contact metric manifolds, almost product manifolds, almost
para-contact metric manifolds, Golden Riemannian manifolds and metallic
Riemannian manifolds. We note that throughout this paper all manifolds
and bundles, along with sections and connections, are assumed to be of
class $C^{\infty }$. A map is always a $C^{\infty }$ map between manifolds.

Let $M'$ be a $2n$-dimensional real manifold. An almost complex structure
$J'$ on $M'$ is a tensor field $J':TM'\rightarrow TM'$ such that
\begin{equation}
J'^2=-I, \label{acm1}
\end{equation}
where $I$ is the identity transformation. Then $(M',J')$ is called almost
complex manifold \cite{YK}.

A smooth map $\varphi:M'_1\rightarrow M'_2$ between almost complex
manifolds $(M'_1,J'_1)$ and $(M'_2,J'_2)$ is called an almost complex
(or holomorphic) map if $d\varphi(J'_1X)=J'_2d\varphi(X)$ for $X\in\Gamma(TM'_1)$,
where $J'_1$ and $J'_2$ are complex structures of $M'_1$ and $M'_2,$
respectively.\\

An $(2n+1)-$dimensional differentiable manifold $M$ is said to have an almost
contact structure $(\phi ,\xi ,\eta )$ if it carries a tensor field $\phi $
of type $(1,1)$, a vector field $\xi $ and $1-$ form $\eta $ on $M$
respectively such that%
\begin{equation}
\phi ^{2}=-I+\eta \otimes\xi,\ \phi\xi =0,\ \eta\circ\phi =0,\ \eta(\xi) =1,  \label{acmm1}
\end{equation}%
where $I$ is the identity transformation. The almost contact structure is
said to be normal if $N+d\eta \otimes \xi =0,$ where $N$ is the Nijenhuis
tensor of $\phi .$ Suppose that a Riemannian metric tensor $g$ is given in $%
M $ and satisfies the condition%
\begin{equation}
g(\phi X,\phi Y)=g( X,Y)-\eta(X)\eta(Y),\ \ \eta(X)=g(X,\xi).
\label{acmm2}
\end{equation}
Then $(M,\phi,\xi,\eta,g)$ is called an almost contact
metric manifold \cite{B1}.\\

Let $N$ be an $n-$dimensional manifold with a tensor of type $(1,1)$ such that
\begin{equation}
\varphi^{2}=I,  \label{apm1}
\end{equation}%
where $I$ is the identity transformation. Then we say that $N$ is an almost
product manifold with almost product structure $\varphi.$ We put%
\begin{equation}
Q=\frac{1}{2}(I+\varphi),\ \ \ Q^{^{\prime}}=\frac{1}{2}(I-\varphi).  \label{apm2}
\end{equation}%
Then we have%
\begin{equation}
Q+Q^{^{\prime }}=I,\ Q^{2}=Q,\ Q^{^{\prime}2}=Q^{^{\prime}},\ QQ^{^{\prime}}=Q^{\prime }Q=0  \label{apm3}
\end{equation}%
and%
\begin{equation}
\varphi=Q-Q^{\prime }.  \label{apm4}
\end{equation}%
If an almost product manifold $N$ admits a Riemannian metric $g$ such that%
$$
g(\varphi X,\varphi Y)=g(X,Y)
$$%
for any vector fields $X$ and $Y$ on $N,$ then $N$ is called an almost
product Riemannian manifold \cite{YK}.\\

Let $N^{\prime}$ be a $(2n + 1)-$dimensional smooth manifold, $\phi^{\prime}$ a (1, 1)-tensor field called the structure endomorphism,
$\xi^{\prime}$ a vector field called the characteristic vector field, $\eta^{\prime}$ a 1 -form called the paracontact form, and $g^{\prime}$ a pseudo Riemannian
metric on $N^{\prime}$ of signature $(n + 1, n)$. In this case, we say that $(\phi^{\prime}, \xi^{\prime}, \eta^{\prime}, g^{\prime})$ defines an almost
paracontact metric structure on $N^{\prime}$ if:
\begin{equation}
\phi^{\prime2}=I-\eta^{\prime}\otimes\xi^{\prime},\ \eta
^{\prime}(\xi^{\prime})=1,\ g^{\prime}(\phi^{\prime}X,\phi^{\prime}Y)=g^{\prime
}(X,Y)-\eta^{\prime}(X)\eta^{\prime}(Y). \label{apcmm1}
\end{equation}
From the definition it follows $\phi^{\prime}(\xi^{\prime})=0,\,\,\eta^{\prime}\circ\phi^{\prime}=0,\,\,\eta^{\prime}(X)=g^{\prime}(X,\xi^{\prime}),\,\,g^{\prime}(\xi^{\prime},\xi^{\prime})=1$ and the fact that $\phi^{\prime}$ is $g^{\prime}-$skewsymmetric: $g^{\prime}(\phi^{\prime}X,Y)=−g^{\prime}(\phi^{\prime}Y,X)$. The associated 2-form $\omega(X,Y)=g^{\prime}(X,\phi^{\prime}Y)$ is skew-symmetric and is
called the fundamental form of the almost metric paracontact manifold $(N^{\prime},\phi^{\prime},\xi^{\prime},\eta^{\prime},g^{\prime})$\cite{CL,C,KW,SC,Z}.\\

Let $(\bar{M},g)$ be a Riemannian manifold with a tensor of
type $(1,1)$ such that
\begin{equation}
P^{2}=P+I,  \label{gm1}
\end{equation}%
where $I$ is the identity transformation. We say that the metric $g$ is $P$
compatible if the equality%
\begin{equation}
g(PX,Y)=g(X,PY)  \label{gm2}
\end{equation}%
for all $X,Y\in\Gamma(TM).$ If we substitute $PX$ into $X$ in
(\ref{gm2}) the equation (\ref{gm2}) may also written as%
$$
g(PX,Y)=g(P^{2}X,Y)=g((P+I)X,Y)=g(PX,Y)+g(X,Y).
$$%
The Riemannian metric (\ref{gm2}) is called $P-$compatible and $(\bar{M},P,g)$
is named a Golden Riemannian manifold \cite{CH2}.\\

Let $(\tilde{M},g)$ be a Riemannian manifold with a tensor of
type $(1,1)$ such that
\begin{equation}
J^{2}=pJ+qI,  \label{mm1}
\end{equation}%
where $p,q$ are positive integers and $I$ is the identity operator on the Lie algebra
$\chi(M)$ of the vector fields on $\tilde{M}$. We say that the metric $g$ is $J-compatible$ if:
\begin{equation}
g(JX,Y)=g(X,JY)  \label{mm2}
\end{equation}%
for every $X,Y\in\chi(M),$ which means that $J$ is a self-adjoint operator with respect to
$g.$ This condition is equivalent in our framework with:

$$
g(JX,JY)=p.g(X,JY)+q.g(X,Y).
$$%
The Riemannian metric (\ref{mm2}) is called $J-$compatible and $(\tilde{M},J,g)$
is named a metallic Riemannian manifold \cite{CH1}.\\

Let $(M,g_{M})$ and $(N,g_{N})$ be Riemannian manifolds and suppose that $%
\pi:M \longrightarrow N$ is a smooth map between them. Then the
differential of $\pi_{*}$ of $\pi$ can be viewed a section of the
bundle $Hom(TM,\pi^{-1}TN) \longrightarrow M$, where $\pi^{-1}TN$ is
the pullback bundle which has fibres $(\pi^{-1}TN)_{p}=T_{\pi(p)}N$,
$p\in M$. $Hom(TM,\pi^{-1}TN)$ has a connection $\nabla$ induced from
the Levi-Civita connection $\nabla^{M}$ and the pullback connection. Then
the second fundamental form of $\pi$ is given by
\begin{equation}
(\nabla\pi_{*})(X,Y)=\nabla^{\pi}_{X}\pi_{*}(Y)-\pi_{*}(%
\nabla^{M}_{X}Y) \label{npixy}
\end{equation}
for $X,Y\in\Gamma(TM)$, where $\nabla^{\pi}$ is the pullback connection.
It is known that the second fundamental form is symmetric. A smooth map $%
\pi:(M,g_{M}) \longrightarrow (N,g_{N})$ is said to be harmonic if $%
trace(\nabla\pi_{*})=0$. On the other hand, the tension field of $\pi$
is the section $\tau(\pi)$ of $\Gamma(\pi^{-1}TN)$ defined by
\begin{equation}
\tau(\pi)=div\pi_{*}=\sum_{i=1}^{m}(\nabla\pi_{*})(e_{i},e_{i}), \label{tau}
\end{equation}
where $\{e_{1},...,e_{m}\}$ is the orthonormal frame on $M$. Then it follows
that $\pi$ is harmonic if and only if $\tau(\pi)=0.$ For more information, see \cite{BW}.

\section{Metallic maps between metallic Riemannian manifolds}

In this section, we give a new notion, namely a metallic map, and give a sufficient condition
 for a map between metallic Riemannian manifolds to be harmonic. We also investigate conditions for a metallic map
to be totally geodesic.\\

\noindent{\bf Definition~3.1.~} Let $F$ be a smooth map from a metallic Riemannian manifold $(M_1,g_1,J_1)$
to a metallic Riemannian manifold $(M_2,g_2,J_2).$ Then $F$ is called a metallic map if and only if
following condition is satisfied.
\begin{equation}
F_*J_1=J_2F_*. \label{memap1}
\end{equation}

We provide the following elementary example.\\

\noindent{\bf Example~3.2.~}Consider the following map defined by
\begin{eqnarray*}
  F:R^{4} &\rightarrow& R^{2}\\
  (x_1,x_1,x_3,x_4)&\rightarrow&(\frac{x_1+x_2}{2},\frac{x_3+x_4}{2}).
\end{eqnarray*}
Then, the kernel of $F_*$ is
$$
\mathcal{V}=KerF_*=Span\{Z_1=\frac{\partial}{\partial \,
x_1}-\frac{\partial}{\partial \,
x_2}\,,Z_2=\frac{\partial}{\partial \,
x_3}-\frac{\partial}{\partial \, x_4}\}
$$
and the horizontal distribution is spanned by
$$
\mathcal{H}=(KerF_*)^{\perp}=Span\{H_1=\frac{\partial}{\partial \,
x_1}+\frac{\partial}{\partial \, x_2}\,,H_2=\frac{\partial}{\partial
\, x_3}+\frac{\partial}{\partial \, x_4}\}.
$$
Then considering metallic structures on $R^4$ and $R^2$ defined by
$$
J_1(a_1,a_2,a_3,a_4)=(\sigma_{p,q}a_1,\sigma_{p,q}a_2,(p-\sigma_{p,q})a_3,(p-\sigma_{p,q})a_4)
$$
and
$$
J_2(a_1,a_2)=(\sigma_{p,q}a_1,(p-\sigma_{p,q})a_2),
$$
where $\sigma_{p,q}$ and $p-\sigma_{p,q}$ are eigenvalues of metallic structures \cite{CH1}. Also by direct computations
$F_*(J_1H_1)=J_2F_*(H_1)=(\sigma_{p,q},0)$ and $F_*(J_1H_2)=J_2F_*(H_2)=(0,p-\sigma_{p,q}).$ Thus $F$ is a metallic map.\\

From now on, when we mention a metallic Riemannian manifold, we will assume that its almost metallic structure
is integrable. It means that $\nabla J=0.$\\

We now give a necessary and sufficient condition for a map $F$ to be totally geodesic. We recall that
a map $F$ is totally geodesic if $\nabla F_*=0.$ A geometric interpretation of a totally geodesic map
is that it maps every geodesic in the total manifold into a geodesic in the base manifold in
proportion to arc lengths.\\

\noindent{\bf Theorem~3.3.~} Let $F$ be a metallic map from a metallic Riemannian manifold $(M_1,g_1,J_1)$
to a metallic Riemannian manifold $(M_2,g_2,J_2)$. Then $F$ is totally geodesic if and only if
\begin{equation}
(\nabla F_*)(X,Y)=\frac{1}{q_1}\{(p_2-p_1)J_2(\nabla F_*)(X,Y)+q_2(\nabla F_*)(X,Y)\} \label{memap2}
\end{equation}
for $X,Y\in\Gamma(TM_1).$

\noindent{\bf Proof.} For $X,Y\in\Gamma(TM_1),$ (\ref{npixy}) and (\ref{mm1}) we have
\begin{eqnarray}
(\nabla F_*)(X,Y)=\frac{1}{q_1}\{\nabla^{^F}_{X}F_{*}(J^{2}_{1}Y-p_1J_1Y)-F_{*}(\nabla_{X}J^{2}_1Y-p_1J_1Y)\}.\nonumber
\end{eqnarray}
Then using (\ref{memap1}) we get
\begin{eqnarray}
(\nabla F_*)(X,Y)&=&\frac{1}{q_1}\{\nabla^{^F}_{X}J_{2}F_{*}(J_1Y)-p_1\nabla^{^F}_{X}F_{*}(J_1Y)
-J_{2}F_{*}(\nabla_{X}J_1Y)\nonumber\\
&+&p_{1}F_{*}(\nabla_{X}J_1Y)\}.\nonumber
\end{eqnarray}
Since $J_1$ integrable, using (\ref{memap1}) we obtain
\begin{eqnarray}
(\nabla F_*)(X,Y)&=&\frac{1}{q_1}\{\nabla^{^F}_{X}J^{2}_{2}F_{*}Y-p_1\nabla^{^F}_{X}J_{2}F_{*}Y
-J^{2}_{2}F_{*}(\nabla_{X}Y)\nonumber\\
&+&p_{1}J_{2}F_{*}(\nabla_{X}Y)\}.\nonumber
\end{eqnarray}
Then using (\ref{mm1}), we get
\begin{eqnarray}
(\nabla F_*)(X,Y)&=&\frac{1}{q_1}\{(p_2-p_1)\nabla^{^F}_{X}J_{2}F_{*}Y+q_{2}\nabla^{^F}_{X}F_{*}Y
-(p_2-p_1)J_2F_{*}(\nabla_{X}Y)\nonumber\\
&-&q_2F_{*}(\nabla_{X}Y)\}.\nonumber
\end{eqnarray}
Since $J_2$ integrable, from (\ref{npixy}) we obtain (\ref{memap2}).\\

We now give a sufficient condition for a map between metallic Riemannian manifold to be harmonic.\\

\noindent{\bf Theorem~3.4.~}  Let $F$ be a metallic map from a metallic Riemannian manifold $(M_1,g_1,J_1)$
to a metallic Riemannian manifold $(M_2,g_2,J_2)$. Then $F$ is a harmonic map if $q_1-q_2\neq (p_2-p_1)\sigma {(p_2,q_2)}.$\\

\noindent{\bf Proof.} Let $\{e_1,e_2,...,e_m\}$ be a basis of $T_{x}M_1,\,\,x\in M_1.$ Then from (\ref{tau}) and (\ref{memap2}),
we have
\begin{eqnarray}
q_1\sum_{i=1}^{m}(\nabla F_{*})(e_{i},e_{i})&=&(p_2-p_1)J_2\sum_{i=1}^{m}(\nabla F_*)(e_{i},e_{i})+q_2\sum_{i=1}^{m}(\nabla F_*)(e_{i},e_{i})\nonumber\\
q_1\tau(F)&=&(p_2-p_1)J_2\tau(F)+q_2\tau(F)\nonumber \\
(q_1-q_2)\tau(F)&=&(p_2-p_1)J_2\tau(F).\label{memap3}
\end{eqnarray}
Applying $J_2$ to (\ref{memap3}), we obtain
\begin{equation}
(q_1-q_2)J_2\tau(F)=(p_2-p_1)J^{2}_2\tau(F).\label{memap4}
\end{equation}
Using (\ref{mm1}), we get
\begin{equation}
(q_1-q_2)J_2\tau(F)=(p_2-p_1)p_2J_2\tau(F)+(p_2-p_1)q_2\tau(F).\label{memap5}
\end{equation}
From above equation, we obtain
\begin{equation}
\{(q_1-q_2)-p_2(p_2-p_1)\}J_2\tau(F)=(p_2-p_1)q_2\tau(F).\label{memap6}
\end{equation}
Now, taking inner product with $(p_2-p_1)q_2$ in (\ref{memap4}) we get
\begin{equation}
(p_2-p_1)(q_1-q_2)q_2\tau(F)=(p_2-p_1)^2q_2J_2\tau(F). \label{memap7}
\end{equation}
In a similar way, taking inner product with $(q_1-q_2)$ in (\ref{memap6}) we obtain
\begin{equation}
(q_1-q_2)\{(q_1-q_2)-p_2(p_2-p_1)\}J_2\tau(F)=(q_1-q_2)(p_2-p_1)q_2\tau(F).\label{memap8}
\end{equation}
Substracting (\ref{memap7}) and (\ref{memap8}), we get
\begin{equation}
\{(p_2-p_1)^2q_2-(q_1-q_2)^2+p_2(p_2-p_1)(q_1-q_2)\}J_2\tau(F)=0. \label{memap9}
\end{equation}
From (\ref{memap9}), since $J_2$ is a isomorphism, $F$ is harmonic map if $(q_1-q_2)\neq (p_2-p_1)\sigma_{(p_2,q_2)}$. This proof is complete.\\

One can see that contrary to the Golden case, metallic map is not always harmonic.\\

\noindent{\bf Remark~3.5.~} We note that for any $C^2$ real valued function $f$ defined
on an open subset of a Riemannian manifold $M,$ the equation $\bigtriangleup f=0$ is called
Laplace's equation and solutions are called harmonic functions on $U.$ Let $F:M\rightarrow N$
be a smooth map between Riemannian manifolds. Then $F$ is called a harmonic morphism if, for
every harmonic function $f:V\rightarrow R$ defined an open subset $V$ of $N$ with $F^{-1}(V)$
non-empty, the composition $f\circ F$ is harmonic on $F^{-1}(V).$ A smooth map $F:M\rightarrow N$
between Riemannian manifolds is harmonic morphism if and only if $F$ is both harmonic and horizontally
weakly conformal \cite{ES}, \cite{F} and \cite{I}. In this respect, from above theorem a metallic map is candidate for a
harmonic morphism.\\

\section{Certain constancy conditions for maps between metallic Riemannian manifolds and manifolds equipped with other differential structures}

In this section we investigate constancy of certain maps between metallic Riemannian
manifolds and manifolds equipped with other differential structures
by imposing holomorphic-like conditions. We first check the situation for a map between
metallic Riemannian manifolds and Golden Riemannian manifolds.\\

\noindent{\bf Theorem~4.1.~} Let $F$ be a smooth map from a metallic Riemannian manifold $(M_1,g_1,J_1)$
to a Golden Riemannian manifold $(\bar{M},P,g)$ such that $F_*J_1=PF_*$ is satisfied. Then $F$ is
a constant map if $p\neq q(\phi)+(1-\phi)$ and $p\neq q(1-\phi)+\phi.$\\

\noindent{\bf Proof.} Let $(M_1,g_1,J_1)$ be a metallic Riemannian manifold and $(\bar{M},P,g)$
a Golden Riemannian manifold. Suppose that $F:M_1\rightarrow \bar{M}$ satisfies
\begin{equation}
F_*(J_1X)=PF_*(X),\,\,\,X\in\Gamma(TM_1). \label{mg1}
\end{equation}
Then apply $P$ to the above equation and using (\ref{gm1}) and (\ref{mm1}), we have
\begin{equation}
pF_*(J_1X)+qF_*(X)=PF_*(X)+F_*(X),\,\,\,X\in\Gamma(TM_1). \label{mg2}
\end{equation}
Using (\ref{mg1}), we get
\begin{equation}
(p-1)F_*(J_1X)=(1-q)F_*(X),\,\,\,X\in\Gamma(TM_1). \label{mg3}
\end{equation}
Applying $P$ to (\ref{mg3}) again and using (\ref{mg1}), we have
\begin{equation}
p^2F_*(J_1X)+pqF_*(X)+qF_*(J_1X)=pF_*(J_1X)+qF_*(X)+F_*(J_1X) \label{mg4}
\end{equation}
for $X\in\Gamma(TM_1).$ From above equation, we get
\begin{equation}
(p^2+q-p-1)F_*(J_1X)=(q-pq)F_*(X),\,\,\,X\in\Gamma(TM_1). \label{mg5}
\end{equation}
From (\ref{mg3}) and (\ref{mg5}) we obtain
\begin{equation}
(p^2-q^2+3q-pq-p-1)F_*(X)=0,\,\,\,X\in\Gamma(TM_1). \label{mg6}
\end{equation}
From (\ref{mg6}), $F$ is a constant map if $p\neq q(\phi)+(1-\phi)$ and $p\neq q(1-\phi)+\phi.$\\

In a similar way, we have the following result.\\

\noindent{\bf Theorem~4.2.~} Let $F$ be a smooth map from a Golden Riemannian manifold $(\bar{M},P,g)$
to a metallic Riemannian manifold $(M_1,g_1,J_1)$ such that $F_*P=J_1F_*$ is satisfied. Then $F$ is
a constant map if $p\neq q(\phi)+(1-\phi)$ and $p\neq q(1-\phi)+\phi.$\\

We now check a similar situation for a map between metallic Riemannian manifolds and almost product
manifolds.\\

\noindent{\bf Theorem~4.3.~} Let $F$ be a smooth map from a metallic Riemannian manifold $(M_1,g_1,J_1)$
to an almost product manifold $(N,g,\varphi)$ such that $F_*J_1=\varphi F_*$ is satisfied. Then $F$ is
a constant map if $p\neq\mp \sqrt{q-1}.$\\

\noindent{\bf Proof.} Let $(M_1,g_1,J_1)$ be a metallic Riemannian manifold and $(N,g,\varphi)$
an almost product manifold. Suppose that $F:M_1\rightarrow N$ satisfies
\begin{equation}
F_*(J_1X)=\varphi F_*(X),\,\,\,X\in\Gamma(TM_1). \label{ma1}
\end{equation}
Then apply $\varphi$ to the above equation and using (\ref{apm1}) and (\ref{mm1}), we have
\begin{equation}
pF_*(J_1X)+qF_*(X)=F_*(X),\,\,\,X\in\Gamma(TM_1). \label{ma2}
\end{equation}
From above equation, we get
\begin{equation}
pF_*(J_1X)=(1-q)F_*(X),\,\,\,X\in\Gamma(TM_1). \label{ma3}
\end{equation}
Applying $\varphi$ to (\ref{ma2}) again and using (\ref{ma1}), we have
\begin{equation}
p^2F_*(J_1X)+pqF_*(X)+qF_*(J_1X)=F_*(J_1X),\,\,\, X\in\Gamma(TM_1). \label{ma4}
\end{equation}
From above equation, we get
\begin{equation}
(p^2+q-1)F_*(J_1X)=-pqF_*(X),\,\,\,X\in\Gamma(TM_1). \label{ma5}
\end{equation}
From (\ref{ma3}) and (\ref{ma5}) we obtain
\begin{equation}
(p^2-q^2+2q-1)F_*(X)=0,\,\,\,X\in\Gamma(TM_1). \label{ma6}
\end{equation}
From (\ref{ma6}), $F$ is a constant map if $p\neq\mp \sqrt{q-1}.$\\

In a similar way, we have the following result.\\

\noindent{\bf Theorem~4.4.~} Let $F$ be a smooth map from an almost product manifold $(N,g,\varphi)$
to a metallic Riemannian manifold $(M_1,g_1,J_1)$ such that $F_*\varphi=J_1F_*$ is satisfied. Then $F$ is
a constant map if and only if $p\neq\mp \sqrt{q-1}.$\\

\noindent{\bf Remark~4.5.~} The equality condition of Theorem 4.4 is satisfies for copper case due to $p=1$ and $q=2.$
So there may be non-constants map between metallic Riemannian manifold and almost product manifold. This situation is
quite different from Golden case.\\

We check a similar situation for a map between metallic Riemannian manifolds and almost complex manifolds.\\

\noindent{\bf Theorem~4.6.~} Let $F$ be a smooth map from a metallic Riemannian manifold $(M_1,g_1,J_1)$
to an almost complex manifold $(M',J')$ such that $F_*J_1=J'F_*$ is satisfied. Then $F$ is
a constant map.\\

\noindent{\bf Proof.} Let $(M_1,g_1,J_1)$ be a metallic Riemannian manifold and $(M',J')$
an almost complex manifold. Suppose that $F:M_1\rightarrow M'$ satisfies
\begin{equation}
F_*(J_1X)=J'F_*(X),\,\,\,X\in\Gamma(TM_1). \label{mac1}
\end{equation}
Then apply $J'$ to the above equation and using (\ref{acm1}) and (\ref{mm1}), we have
\begin{equation}
pF_*(J_1X)+qF_*(X)=-F_*(X),\,\,\,X\in\Gamma(TM_1). \label{mac2}
\end{equation}
From above equation, we get
\begin{equation}
pF_*(J_1X)=-(q+1)F_*(X),\,\,\,X\in\Gamma(TM_1). \label{mac3}
\end{equation}
Applying $J'$ to (\ref{mac2}) again and using (\ref{mac1}), we have
\begin{equation}
p^2F_*(J_1X)+pqF_*(X)=-(q+1)F_*(J_1X),\,\,\, X\in\Gamma(TM_1). \label{mac4}
\end{equation}
From above equation, we get
\begin{equation}
(p^2+q+1)F_*(J_1X)=-pqF_*(X),\,\,\,X\in\Gamma(TM_1). \label{mac5}
\end{equation}
From (\ref{mac3}) and (\ref{mac5}) we obtain
\begin{equation}
(p^2+(q+1)^2)F_*(X)=0,\,\,\,X\in\Gamma(TM_1). \label{mac6}
\end{equation}
which shows that $F$ is constant.\\

In a similar way, we have the following result.\\

\noindent{\bf Theorem~4.7.~} Let $F$ be a smooth map from an almost complex manifold $(M',J')$
to a metallic Riemannian manifold $(M_1,g_1,J_1)$ such that $F_*J'=J_1F_*$ is satisfied. Then $F$ is
a constant map.\\

The following result shows that a smooth map satisfying a compatible condition between metallic
Riemannian manifolds and almost contact metric manifolds is also constant.\\

\noindent{\bf Theorem~4.8.~} Let $F$ be a smooth map from a metallic Riemannian manifold $(M_1,g_1,J_1)$
to an almost contact metric manifold $(M,\phi,\xi,\eta,g)$ such that $F_*J_1=\phi F_*$ is satisfied. Then $F$ is
a constant map.\\

\noindent{\bf Proof.} Let $(M_1,g_1,J_1)$ be a metallic Riemannian manifold and $(M,\phi,\xi,\eta,g)$
an almost contact metric manifold. Suppose that $F:M_1\rightarrow M$ satisfies

\begin{equation}
F_*(J_1X)=\phi F_*(X),\,\,\,X\in\Gamma(TM_1). \label{macm1}
\end{equation}
Then apply $\phi$ to the above equation and using (\ref{acmm1}) and (\ref{mm1}), we have
\begin{equation}
pF_*(J_1X)+(q+1)F_*(X)=\eta(F_*(X))\xi,\,\,\,X\in\Gamma(TM_1). \label{macm2}
\end{equation}
Then applying $\phi$ to (\ref{macm2}) again and using (\ref{macm1}) and (\ref{acmm1}), we have
\begin{equation}
(p^2+q+1)F_*(J_1X)=-pqF_*(X),\,\,\,X\in\Gamma(TM_1). \label{macm3}
\end{equation}
From (\ref{macm2}) and (\ref{macm3}), we obtain
\begin{equation}
F_*(X)\left(\frac{p^2+(q+1)^2}{p^2+q+1}\right)=\eta(F_*(X))\xi,\,\,\, X\in\Gamma(TM_1). \label{macm4}
\end{equation}
Again applying $\phi$ to (\ref{macm4}), we get
\begin{equation}
\phi F_*(X)\left(\frac{p^2+(q+1)^2}{p^2+q+1}\right)=0,\,\,\,X\in\Gamma(TM_1). \label{macm5}
\end{equation}
Then applying $\phi$ to (\ref{macm5}), we get
\begin{equation}
-F_*(X)\left(\frac{p^2+(q+1)^2}{p^2+q+1}\right)+\eta\left(F_*(X)\frac{p^2+(q+1)^2}{p^2+q+1}\right)\xi=0,\,\,\,X\in\Gamma(TM_1). \label{macm6}
\end{equation}
Again applying $\phi$ to (\ref{macm6}), we get
\begin{equation}
F_*(X)\left(\frac{p^2+(q+1)^2}{p^2+q+1}\right)=0,\,\,\,X\in\Gamma(TM_1). \label{macm7}
\end{equation}
From (\ref{macm7}), $F$ is a constant map.\\

In a similar way, we have the following result.\\

\noindent{\bf Theorem~4.9.~} Let $F$ be a smooth map from an almost contact metric manifold $(M,\phi,\xi,\eta,g)$
to a metallic Riemannian manifold $(M_1,g_1,J_1)$ such that $F_*\phi=J_1F_*$ is satisfied. Then $F$ is
a constant map.\\

Finally, we check the same problem for almost para-contact metric manifolds.\\

\noindent{\bf Theorem~4.11.~} Let $F$ be a smooth map from a metallic Riemannian manifold $(M_1,g_1,J_1)$
to an almost para-contact metric manifold $(N^{\prime},\phi^{\prime},\xi^{\prime},\eta^{\prime},g^{\prime})$
such that $F_*J_1=\phi^{\prime}F_*$ is satisfied. Then $F$ is a constant map if $p^2\neq(q-1)^2.$\\

\noindent{\bf Proof.} Let $(M_1,g_1,J_1)$ be a metallic Riemannian manifold and
$(N^{\prime},\phi^{\prime},\xi^{\prime},\eta^{\prime},g^{\prime})$ an almost para-contact metric manifold.
Suppose that $F:M_1\rightarrow N'$ satisfies

\begin{equation}
F_*(J_1X)=\phi^{\prime}F_*(X),\,\,\,X\in\Gamma(TM_1). \label{mapm1}
\end{equation}
Then apply $\phi^{\prime}$ to the above equation and using (\ref{apcmm1}) and (\ref{mm1}), we have
\begin{equation}
pF_*(J_1X)+(q-1)F_*(X)=-\eta(F_*(X))\xi,\,\,\,X\in\Gamma(TM_1). \label{mapm2}
\end{equation}
Then applying $\phi^{\prime}$ to (\ref{mapm2}) again and using (\ref{mapm1}) and (\ref{apcmm1}), we have
\begin{equation}
(p^2+q-1)F_*(J_1X)=-pqF_*(X),\,\,\,X\in\Gamma(TM_1). \label{mapm3}
\end{equation}
From (\ref{mapm2}) and (\ref{mapm3}), we obtain
\begin{equation}
F_*(X)\left(\frac{(q-1)^2-p^2}{p^2+q-1}\right)=-\eta(F_*(X))\xi,\,\,\, X\in\Gamma(TM_1). \label{mapm4}
\end{equation}
Again applying $\phi^{\prime}$ to (\ref{mapm4}), we get
\begin{equation}
\phi^{\prime}F_*(X)\left(\frac{(q-1)^2-p^2}{p^2+q-1}\right)=0,\,\,\,X\in\Gamma(TM_1). \label{mapm5}
\end{equation}
Then applying $\phi^{\prime}$ to (\ref{mapm5}), we get
\begin{equation}
F_*(X)\left(\frac{(q-1)^2-p^2}{p^2+q-1}\right)-\eta\left(F_*(X)\frac{(q-1)^2-p^2}{p^2+q-1}\right)\xi=0 \label{mapm6}
\end{equation}
for $X\in\Gamma(TM_1).$ Again applying $\phi^{\prime}$ to (\ref{mapm6}), we get
\begin{equation}
F_*(X)\left(\frac{(q-1)^2-p^2}{p^2+q-1}\right)=0,\,\,\,X\in\Gamma(TM_1). \label{mapm7}
\end{equation}
From (\ref{mapm7}), $F$ is a constant map if $p^2\neq(q-1)^2.$\\

In a similar way, we have the following result.\\

\noindent{\bf Theorem~4.10.~} Let $F$ be a smooth map from an almost para-contact metric manifold $(N^{\prime},\phi^{\prime},\xi^{\prime},\eta^{\prime},g^{\prime})$
to a metallic Riemannian manifold $(M_1,g_1,J_1)$ such that $F_*\phi^{\prime}=J_1F_*$ is satisfied.
Then $F$ is a constant map if $p^2\neq(q-1)^2.$\\

\noindent{\bf Remark~4.11.~} The equality condition of Theorem 4.9 is possible for copper case due to $p=1$ and $q=2.$
Thus it is possible to find non-constant map between such manifolds.\\


\begin{thebibliography}{99}
\bibitem{BA}B. \d{S}ahin and M. A. Akyol,  {\it Golden maps between Golden Riemannian manifolds and constancy of certain maps,} Math. Commun. 19(2014), 333-342.
\bibitem{BW}P. Baird and J. C. Wood, {\it Harmonic Morphisms Between Riemannian Manifolds,} London Mathematical Society Monographs, No. 29, Oxford
University Press, The Clarendon Press, Oxford, 2003.

\bibitem{B1}D. E. Blair, {\it Riemannian Geometry of Contact and Symplectic Manifolds,} Birkhauser, Boston, (2002).

\bibitem{CH1} M. Crasmareanu and C-E. Hret\c{c}anu, {\it Metallic structures on Riemannian manifolds,} Revista de la uni\'{o}n
 matem\'{a}tica argentina vol:54, No:2, (2013), 15-27.

\bibitem{CH2}M. Crasmareanu and C-E. Hret\c{c}anu, {\it Golden differansiyel geometry,} Chaos, Solitons \& Fractals volume 38, issue 5, (2008), 1229-1238.

\bibitem{CH3}M. Crasmareanu and C-E. Hret\c{c}anu, {\it Applications of the Golden Ratio on Riemannian Manifolds,} Turkish J. Math. 33, no. 2, (2009), 179-191.

\bibitem{CL}M. Crasmareanu and L. I. Pi\c{s}coran, {\it Invariant distributions and holomorphic vector fields in paracontact geometry,} Turkish J. Math. 39, no. 4,
(2015), 467-476.

\bibitem{C} H. \c{C}ay\i r, {\it Some Notes on Lifts of Almost Paracontact Structures,} American Review of Mathematics and Statistics. 3(1), (2015), 52-60.

\bibitem{E} E. Sadettin, {\it On product and Golden structures and harmonicty}, arXiv:1612.07731v1 [math.DG].

\bibitem{ES}J. Eells and H. J. Sampson, {\it Harmonic mappings of Riemannian Manifolds,} Amer. J. Math. 86 (1964), 109-160.

\bibitem{F}B. Fuglede, {\it Harmonic morphisms between Riemannian Manifolds,} Ann. Inst. Fourier(Grenoble). 28 (1978), 107-144.

\bibitem{I}T. Ishihara, {\it A mapping of Riemannian manifolds which preserves harmonic functions,} J. Math. Kyoto Univ. 19 (1979), 215-229.

\bibitem{KW}S. Kaneyuki and F. L. Willams, {\it Almost paracontact and parahodge structures on manifolds,} Nagoya Math. J. 99 (1985), 173-187.

\bibitem{GCS}A. Gezer, N. Cengiz, A. Salimov, {\it On integrability of Golden Riemannian structures,} Turk J. Math. 37 (2013), 693-703.

\bibitem{OCT}M. \"{O}zkan, A. \c{C}\i tlak, E. Taylan, {\it Prolongations of Golden Structure to Tangent Bundle of Order 2}, Gazi University Journal of Science, 28(2): 253-258, (2015).

\bibitem{SC} A. Salimov, H. \c{C}ay\i r, {\it Some Notes On Almost Paracontact Structures,} Comptes Rendus de l'Academie Bulgare des Sciences. 66(3), (2013), 331-338.

\bibitem{S1}V.W. de Spinadel, {\it On characterization of the onset to chaos,} Chaos, Solitions \& Fractals, 8(10) (1997), 1631-1643.

\bibitem{S2}V.W. de Spinadel, {\it The metallic means family and multifractal spectra,} Nonlinear Anal. Ser. B: Real World Appl.
36(6) (1999), 721-745.

\bibitem{S3}V.W. de Spinadel, {\it The family of metallic means,} Vis. Math. 1, 3(1999), http.//vismath1.tripod.com/spinadel/.

\bibitem{S4}V.W. de Spinadel, {\it The metallic means family and forbidden symmetries,} Int. Math. J. 2(3) (2002) 279-288.

\bibitem{S5}V.W. de Spinadel, {\it The metallic means family and renormalization group techniques,} Proc. Steklov Inst. Math., Control in Dynamic
systems, suppl. 1 (2000) 194-209.

\bibitem{YK}K. Yano and M. Kon, {\it Structures on Manifolds}. World Scientific, Singapuore, Series in pure mathematics-Volume 3, (1984).

\bibitem{Z}S. Zamkovoy, {\it Canonical connections on paracontact manifolds,} Ann. Glob. Anal. Geom. 36 (2009), 37-60.
\end{thebibliography}
\end{document}